\documentclass[12pt]{amsart}

\usepackage{amsmath,amssymb,amscd,amsfonts,verbatim}

\usepackage[mathscr]{eucal}

\newcommand{\pionealg}{\pi_1^{\mathrm{alg}}}
\newtheorem{thm}{Theorem}[section]
\newtheorem{lem}[thm]{Lemma}
\newtheorem{rem}[thm]{Remark}
\newtheorem{prop}[thm]{Proposition}

\newtheorem{assu-nota}[thm]{Assumption--Notation}
\newtheorem{ex}{Example}[section]
\theoremstyle{remark}

\newcommand{\Z}{\mathbb Z}

\newcommand{\pp}{\mathbb P}
\DeclareMathOperator{\Aut}{Aut}

\newcommand{\Si}{\Sigma}

\newcommand{\fie}{\varphi}

\newcommand{\OO}{\mathcal{O}}

\numberwithin{equation}{section}

\title [ The order of finite algebraic fundamental groups... ]
{ The order of finite algebraic fundamental groups of surfaces with $K^2\le3\chi-2$
}
\author{ Margarida Mendes Lopes \and Rita Pardini}

\thanks{The first author is a member of the Center for Mathematical
Analysis, Geometry and Dynamical Systems and the second author is a member
of G.N.S.A.G.A.-I.N.d.A.M. This research was partially supported by the
Italian project ``Geometria sulle variet\`a algebriche" (PRIN COFIN 2004)
and by FCT (Portugal) through program POCTI/FEDER and Project
POCTI/MAT/44068/2002.}

\begin{document}

\begin{abstract} In this note we study the structure of $\pionealg(S)$ for minimal surfaces of general type  $S$ satisfying $K_S^2\leq 3\chi-2$ and not having  any irregular \'etale cover. We show that, if $K_S^2\leq 3\chi-2$, then $|\pionealg(S)|\leq 5$, and equality only occurs if $S$ is a Godeaux surface. We also show that if $K_S^2\leq 3\chi-3$ and $\pionealg(S)\neq \{1\}$, then 
$\pionealg(S)=\Z_2$, or $\pionealg(S)=\Z_2^2$ or  $\pionealg(S)=\Z_3$. Furthermore  in this last case one has: $2\le \chi\le 4$, $K^2=3\chi-3$ and these possibilities do occur.
\medskip

\end{abstract}
\maketitle

\section { Introduction}

In this note we study the structure of $\pionealg(S)$ for minimal surfaces of general type  $S$ satisfying $K_S^2\leq 3\chi-2$ and not having  any irregular \'etale cover.

In \cite{mp} we have shown, among other things, that  if $S$ has no irregular \'etale cover and $K^2_S\leq 3\chi-1$ then the order of $\pionealg(S)$ is less than or equal to 9 and equality is only possible if $\chi(S)=1$. In this note we show  some more results on the structure of  $\pionealg(S)$.

We want to remark that most of the present  results are  somehow  implicit in \cite{xiaohyp}, but we think it is worthwhile spelling them out.

We prove the following:
\begin{thm}\label{main} Let $S$ be a minimal algebraic surface of general type such that $K^2\leq 3\chi-2$ not having any  irregular \'etale cover. Then 
the order of $\pionealg(S)$ is at most 5, and equality only occurs if 
$\chi=1$ and $K^2=1$ (i.e. $S$ is a numerical Godeaux surface).
\end{thm}

\begin{thm}\label{main1} Let $S$ be a minimal algebraic surface of general type such that $K^2\leq 3\chi-3$ not having any  irregular \'etale cover.  If  $\pionealg(S)\ne \{1\}$, then there are the following posibilities:
 \begin{enumerate}
 \item $\pionealg(S)=\Z_2$;
 \item $\pionealg(S)=\Z_2^2$;
\item $\pionealg(S)=\Z_3$. In this case one has: $2\le \chi\le 4$, $K^2=3\chi-3$.\end{enumerate}
\end{thm}

\begin{rem}\label{sharp} \rm We remark that there are examples of surfaces with $\pionealg=\Z_3$ for all the values of $\chi$ and $K^2$ given in (iii) of Theorem  \ref{main1}.

For  $\chi(S)=2$ see \cite{mura}, whilst for $\chi(S)=3$ it is enough to take $S$ as an \'etale triple cover of the Campedelli surfaces with fundamental group of order 9  described in \cite{pi9}.

An example for $\chi=4$ is described in \ref{example}.
\end{rem}

\begin{rem} \rm  For details on Godeaux surfaces see for instance  the introduction and the references in  \cite{ccm}. 
\end{rem}

\noindent {\bf Notation} We work over the complex numbers. All varieties are projective algebraic. All the notation we use is standard in algebraic geometry. 
We just recall the definition of the numerical invariants of a smooth surface $S$:  the self-intersection number $K^2_S$ of the canonical divisor $K_S$, the {\em geometric genus} $p_g(S):=h^0(K_S)=h^2(\OO_S)$, the {\em irregularity} $q(S):=h^0(\Omega^1_S)=h^1(\OO_S)$ and the {\em holomorphic Euler characteristic} $\chi(S):=1+p_g(S)-q(S)$.

\medskip

\noindent {\bf Acknowledgments} The first author wants to thank the organizers of the Workshop ``Algebraic Geometry and Topology'' for the invitation to the workshop and the wonderful hospitality.

\section{Some useful facts}
We will use the following  fundamental facts:

\begin{prop} \label{be}{\rm (\cite[Cor. 5.8]{beauville}, cf. Proposition 4.1 of \cite{mp})} Let $Y$ be a surface of general type such that
the canonical map of $Y$ has degree~$2$ onto a rational surface. If $G$
is a group that acts freely on $Y$, then $G=\Z_2^r$, for some $r$.
\end{prop}

\begin{lem}\label{free} Let $Y$ be a regular surface of general type,  let $G\ne \{1\} $ be  a finite group   that acts freely on $Y$ and  let $|F|$ be a $G$-invariant free pencil $|F|$ of curves of genus $g(F)\leq 4$. Then only the following possibilities can occur:
\begin{enumerate}
\item $G=\Z_2^2$, $g(F)=3$ and $G$ acts faithfully on $|F|$;
\item $G=\Z_3$, $g(F)=4$;
 \item $G=\Z_2$, $g(F)=3$. 
\end{enumerate}
\end{lem} 
\begin{proof}  
Assume that such a pencil $|F|$ exists, let $H$ be the subgroup of $G$ consisting of the elements that act
trivially on $|F|$ and let $h$ be the order of $H$. Set $Y':=Y/H$ and $S:=Y/G$. The pencil $|F|$ induces a free pencil $|F'|$ on $Y'$ with general fibre $F':=F/H$. Denote by $f'\colon Y'\to\pp^1$ the morphism given by $|F'|$. There is a cartesian diagram:
\begin{equation}
\begin{CD} Y'@>>>S\\ @V{f'}VV @VV{f}V\\
\pp^1@>{p}>>\pp^1
\end{CD}
\end{equation}
where $p$ is a $G/H$-cover and the general fibre of $f$ is also equal to $F'$.  The fibres of $f$ over the branch points of $p$ are multiple fibres, of multiplicity equal to the branching order. 
Hence, if $G/H$ is nontrivial, then $f$ has multiple fibres.

Since $g(F')=1+(g(F)-1)/h$ (recall that $H$ acts freely), $g(F)\le 4$ and $S$ is of general type, we get $h\le 3$.

Assume that   $|G|\ge 4$.  Then  $G/H$ is not the trivial group and $f$ has multiple fibres, hence $g(F')>2$ by the adjunction formula. 
So the only possibility is $h=1$,  $g(F)=3$ or $4$ and $G$ is isomorphic to a subgroup of $\Aut |F|=\Aut \pp^1$. 

If $g(F)=3$, then by the adjunction formula the multiple fibres of $f$ are double fibres. Since every automorphism of $\pp^1$ has fixed points, this implies that every element of $G$ has order 2. Hence $G=\Z_2^r$ for some $r$. Since $\Aut \pp^1$ does not contain a subgroup isomorphic to $\Z_2^3$, we have $r\le 2$. 

If $g(F)=4$, then by the adjunction formula the multiple fibres of $f$ are triple fibres.  Since every automorphism of $\pp^1$ has fixed points, this implies that every element of $G$ has order 3. It is well known (cf. \cite{sfera} or \cite{enriques}) that a finite subgroup of $\Aut \pp^1$ is isomorphic to one of the following: $\Z_n$, $\Z_2^2$, the dihedral group $D_n$, the symmetric group $S_4$, the alternating groups $A_4$ and $A_5$. It follows that $G=\Z_3$ in this case.
So we have proven that the only possibility for $|G|\ge 4$ is (i).

Statements (ii) and (iii) now follow from the adjunction formula.
\end{proof}

\section{The proof of Theorem \ref{main}}

In this section we denote by $S$ a surface satisfying the assumptions of  Theorem \ref{main}, i.~e., a surface such that $K_S^2=3\chi-m$, where $m\geq 2$, having no irregular \'etale cover. Note that, by  Theorem 1.3 of \cite{mp},  $|\pionealg(S)|\leq 8$.

We divide the  proof of Theorem \ref{main} in two  steps.

\smallskip
\noindent{\bf Step 1:} {\em  $|\pionealg(S)|\leq 5$.}

\smallskip

Assume by contradiction that $d:=|\pionealg(S)|> 5$  and let
$Y\to S$ be the corresponding \'etale G-cover of degree $d$.

Since $K_Y^2=dK_S^2=3d\chi(S)-dm=3\chi(Y)-dm< 3\chi(Y)-10$, by \cite{beauville} the canonical map of $Y$ is generically finite of degree 2 and its image is a ruled surface $\Si$. The surface $Y$ is regular and so $\Si$ is a rational surface. Since $|\pionealg(S)|\le 8$,  by Proposition \ref{be} we conclude that $G=\Z_2^3$.

So $d=8$ and $K_Y^2=3\chi(Y)-8m$.

Assume that $m\geq 3$. Then  $K_Y^2\leq 3\chi-24< 3(\chi-\frac{43}{8})$ and so,  by Theorem 1.1 of \cite{xiaohyp}, $Y$ has a unique free pencil of hyperelliptic curves of genus $g\leq 3$. Since this pencil is necessarily G-invariant we have a contradiction to Lemma \ref{free}.

Suppose that $m=2$. Then $K_Y^2=3\chi(Y)-16<3(\chi(Y)-5)$. If $\chi(Y)\geq 24$ (i.~e. $\chi(S)\geq 3$),   again by  Theorem 1.1 of \cite{xiaohyp}, $Y$ has a unique free pencil of hyperelliptic curves of genus $g\leq 3$ and again we have the same contradiction.
If $\chi(Y)=16$ (i.~e. $\chi(S)=2$), then, by the Example on page 133 of \cite{xiaohyp}, either $Y$ has a unique pencil of hyperelliptic curves of genus $\leq 3$ or $Y$ is a double cover of $\pp^2$ with branch locus a curve of degree 14 with at most non-essential singularities. The first case can be excluded as before whilst  the second case  is excluded because $Y$ admits no free involution (see the proof of Lemma 3 of \cite{xiaohyp}).

If $\chi(S)=1$ then by Noether's inequality $S$ cannot have an \'etale cover of degree 8 (the resulting $Y$ would satisfy $K_Y^2=8<2\chi-6=10$).

\bigskip
\noindent{\bf Step 2:} {\em If  $|\pionealg(S)|=5$, then $S$ is a Godeaux surface.}

\smallskip

Assume $|\pionealg(S)|= 5$ and let
$Y\to S$ be the corresponding \'etale G-cover.

Then $K_Y^2=5K_S^2=15\chi(S)-5m=3\chi(Y)-5m$. If $m\geq 3$, then $K_Y^2=3\chi(Y)-5m<3\chi(Y)-10$ and so we obtain a contradiction as in Step 1, because $\Z_2^n$ has not order $5$ for any $n$.

If $m=2$, then the regular surface $Y$ satisfies $K_Y^2=3\chi-10=3p_g-7$ and by \cite
{ak} either the canonical map $\fie$ of $Y$ is  birational or it is a generically finite map of degree $2$  on a rational surface. This last case is immediately excluded by Proposition \ref{be}.

So the canonical map $\fie$ of $Y$ is birational. To show the claim we only need  
to show that $\chi(S)\geq 2$ does not occur in this situation.
If $\chi(S)\geq 2$ then $\chi(Y)\geq 10$.  So, by  \cite{ak}, $Y$ again has an unique free pencil of genus 3 curves and therefore,  by  Lemma \ref{free}, $Y$ does not admit a free action by a group of order 5.

\section { The proof of Theorem \ref{main1}}

 Here we prove Theorem \ref{main1}.
We continue assuming that $S$ is a surface having no irregular \'etale cover and now we suppose  that $K^2\leq 3\chi-3$.

By Theorem \ref{main} the order of $\pionealg(S)$ is less than  or equal to 4.
If the order of  $\pionealg(S)$ is $4$, then the corresponding \'etale cover $Y$ satisfies $K_Y^2\leq 3\chi(Y)-12$. Since by \cite{beauville}  the canonical map of $Y$ is 2-1 onto a rational surface, by  Proposition \ref{be} we have $\pionealg(S)=\Z_2^2$.

Assume that $\pionealg(S)=\Z_3$. Then by the same argument  one has  $K_S^2=3\chi(S)-3$, and the corresponding \'etale cover $Y\to S$ satisfies $K_Y^2=3\chi(Y)-9=3p_g(Y)-6$. 

Surfaces with these invariants have been classified by Konno (\cite{konnopg}). If $\chi(Y)\geq 13$ and $q(Y)=0$,  then either the canonical map of $Y$ is generically finite of degree 2 onto a rational surface, or the canonical map is birational and $Y$ has a unique free pencil of genus 3 non-hyperelliptic  curves.
The first possibility does not occur by  Proposition \ref{be} and the second one does not occur by Lemma \ref{free}.
 Therefore,  if the order of $\pionealg(S)$ is 3, then necessarily $\chi(Y)\leq 12$, i.~e. $\chi(S)\leq 4$.
Finally notice that $\chi(S)=1$ does not occur, since by assumption $3\chi(S)-3\ge K^2_S$ and $K^2_S>0$ since $S$ is of general type.
This finishes the proof of Theorem \ref{main1}.

\begin{ex}\label{example}\rm As explained in Remark \ref{sharp}, examples of surfaces with $K^2=3\chi-3$ and $\pionealg=\Z_3$ are known  in the literature for $\chi=2,3$. We describe here an example with $\chi=4$.

Consider homogeneous coordinates $(x_0,x_1, x_2)$ on $\pp^2$ and let $1\in\Z_3$ act on $\pp^2$ by:  $(x_0,x_1,x_2)\mapsto (x_0, \omega x_1, \omega^2 x_2)$, where $\omega\ne 1$ is a cube root of 1. Consider $\pp^{10}$ with homogeneous coordinates $z_1,\dots z_{10}$ and identify $\pp^9$ with the hyperplane of $\pp^{10}$ defined by $z_0=0$. 
Denote by $m_1,\dots m_{10}$ the homogeneous monomials of degree 3 in $x_0, x_1,x_2$. The Veronese embedding of degree 3 $v\colon \pp^2\to \pp^9$ is defined by letting  $z_i=m_i$, $i=1,\dots 10$.  Denote by $\Si$ the image of $v$ and by $K\subset \pp^{10}$ the cone over $\Si$ with vertex the point $P:=(1,0,\dots 0)$. The $\Z_3$-action on $\pp^2$ induces a compatible action on $\Si\subset \pp^9$,  which is diagonal with respect to the coordinates $z_1,\dots z_{10}$. We extend this action to $\pp^{10}$ by letting $1\in \Z_3$ act on $z_0$ as multiplication by $\omega$. Clearly this action fixes the point $P$ and maps the cone $K$ to itself. We claim that the only fixed points of $\Z_3$ on $K$ are $P$ and the image points $Q_0,Q_1, Q_2\in\Si$ of the coordinate points of $\pp^2$. In fact, let $Q\in K$, $Q\ne P$ be a fixed point. Then the line $<P,Q>$ meets $\Si$ in a point fixed by $\Z_3$, namely $Q$ lies on one of the lines $<P,Q_i>$, for some $i\in\{0,1,2\}$. So it is enough to show that $\Z_3$ acts non trivially on the line $<P,Q_i>$ for $i=0,1,2$. This is easy to check, since the $Q_i$ are coordinate points in $\pp^9$ and the only nonzero coordinate of $Q_i$ is the one corresponding to the monomial $x_i^3$, hence $\Z_3$ acts  on it as multiplication by 1. 

Let now $T_1\subset H^0(\OO_{\pp^{10}}(3))$ be the subspace of elements which are fixed by $\Z_3$. Notice that $z_0^3,\dots z_{10}^3$ are elements of $T_1$, hence the system $|T_1|$ is free. Let $V$ be the intersection of $K$ with a general hypersurface of $|T_1|$. Then by Bertini's theorem $V$ is smooth, $\Z_3$-invariant and the $\Z_3$-action on $V$ is free. One has $K_V=\OO_V(1)$, $K^2_V=27$, $p_g(V)=11$, $q(V)=0$ (cf. \cite[\S 4]{konnopg}). The quotient surface $S:=V/\Z_3$ is smooth, minimal of general type with $K^2_S=9$, $q(S)=0$, $\chi(S)=4$. By Theorem \ref{main1}, $V$ is the universal cover of $S$ and $S$ is an example of case (iii) of Theorem \ref{main}. (One can also check directly that $V$ is simply connected).

\end{ex}

\bigskip

\bigskip

\begin{minipage}{13cm}
\parbox[t]{6.7cm}{Margarida Mendes Lopes\\
Departamento de Matem\'atica\\
Instituto Superior T\'ecnico\\
Universidade T{\'e}cnica de Lisboa\\
Av.~Rovisco Pais\\
1049-001 Lisboa, PORTUGAL\\
mmlopes@math.ist.utl.pt
} \hfill
\parbox[t]{5.5cm}{Rita Pardini\\
Dipartimento di Matematica\\
Universit\`a di Pisa\\
Largo B. Pontecorvo, 5\\
56127 Pisa, Italy\\
pardini@dm.unipi.it}
\end{minipage}


\begin{thebibliography}{ABC}
\bibitem[AK]{ak} T. Ashikaga, K. Konno, {\em Algebraic surfaces of general type with $c\sp 2\sb 1=3p\sb g-7$}, Tohoku Math. J. (2) {\bf 42} (1990), no. 4, 517--536. 	
\bibitem[Be]{beauville} A. Beauville, {\em L'application canonique pour les surfaces de type g\'en\'eral}. Inv. Math. {\bf 55} (1979), 121--140.
\bibitem [Bl]{sfera} H.~Blichfeldt, {\em Finite Collineation Groups}, The
University of Chicago Science series, The University of Chicago Press,
Chicago, Illinois, 1917.
\bibitem[CCM]{ccm}
A.~Calabri, C.~Ciliberto, M.~Mendes Lopes,
{\em Numerical Godeaux surfaces with an involution}, to appear in
Trans. A.M.S..
\bibitem[EC]{enriques} F. Enriques, O. Chisini, {\em Teoria geometrica delle
equazioni},
Zanichelli E\-di\-tore, Bologna, 1915.

\bibitem[Ko]{konnopg}
K. Konno, {\em Algebraic surfaces of general type with $c\sp 2\sb 1=3p\sb g-6$},
Math. Ann. {\bf 290} (1991), no. 1, 77--107.
\bibitem[MP1]{pi9} M. Mendes Lopes, R. Pardini, {\em Numerical Campedelli
surfaces with fundamental group of order 9}, preprint  math.AG/0602633.
\bibitem[MP2]{mp} M.Mendes Lopes, R.Pardini, {\em On the algebraic fundamental group of surfaces with $K^2\leq 3\chi$}, preprint math.AG/0512483.
\bibitem[Mur]{mura} M.~Murakami, {\em Minimal algebraic surfaces of general type with $c\sp 2\sb 1=3, p\sb g=1$ and $q=0$, which have non-trivial 3-torsion divisors},  J. Math. Kyoto Univ. {\bf 43}  (2003),  no. 1, 203--215. 
\bibitem[X]{xiaohyp} G. Xiao, {\em Hyperelliptic surfaces of general type
with $K^2<4\chi$}, Manuscripta Math. {\bf 57} (1987), 125--148.
\end{thebibliography}
\end{document}